\documentclass[12pt,usenames,dvipsnames]{article}
\usepackage[a4paper]{anysize}\marginsize{2cm}{2cm}{2cm}{2cm}
\pdfpagewidth=\paperwidth \pdfpageheight=\paperheight
\usepackage{amsfonts,amssymb,amsthm,amsmath,eucal,tabu,hyperref}
\usepackage{xcolor}
\usepackage{float}
\usepackage{abstract}

\usepackage{wrapfig}
\usepackage{pgfplots}
\pgfplotsset{compat=1.15}
\usepackage{mathrsfs}
\usetikzlibrary{arrows}
\usepackage{diagbox}
\usepackage{subcaption}

\usepackage{tikz-cd}

\usepackage[
  backend=biber,
  style=alphabetic,
  giveninits=true,
  maxnames=3,
  minnames=1
]{biblatex}
\theoremstyle{plain}
\newtheorem{thm}{Theorem}[section]

\newtheorem{proposition}[thm]{Proposition}

\theoremstyle{definition}
\newtheorem{definition}[thm]{Definition}
\newtheorem{remark}[thm]{Remark}

\newtheorem{example}[thm]{Example}

\DeclareMathOperator{\Der}{Der}

\DeclareMathOperator{\Jump}{Jump}
\DeclareMathOperator{\vdim}{vdim}

\DeclareMathOperator{\edim}{edim}

\DeclareMathOperator{\Sing}{Sing}

\def\PP{\mathbb{P}}
\def\C{\mathbb{C}}

\def\cA{\mathcal{A}}

\def\cO{\mathcal{O}}
\def\cE{\mathcal{E}}
\def\cL{\mathcal{L}}
\def\cI{\mathcal{I}}

\def\cK{\mathcal{K}}

\def\fm{\mathfrak{m}}
\def\a0n{a_0,\ldots,a_N}
\def\x0n{x_0,\ldots,x_N}
\def\b0n{b_0,\ldots,b_N}
\def\y0n{y_0,\ldots,y_N}
\def\p0n{p_0,\ldots,p_N}

\author{Elena Guardo, Graham Keiper, Grzegorz Malara}
\date{\today}

\begin{document}
\title{Interpolation matrices and jumping lines of logarithmic bundles}
\maketitle

\begin{abstract}
    We study  jumping lines loci of logarithmic bundles associated with finite sets of points in the projective plane. Using the interpolation matrix introduced in \cite{DMTG25}, we describe these loci as the zero sets of explicit determinants depending on parameters $(d,m)$ determined by the number of points. We show that for points in general position the determinant defines an irreducible curve of the expected degree, while for special configurations it acquires fixed components related to the combinatorics of the arrangement. The approach provides a new geometric interpretation of the classical jumping lines of Dolgachev--Kapranov and Barth, and connects them to the framework of unexpected curves and hypersurfaces.


\end{abstract}

\section{Introduction}

The study of the geometry of rank--2 vector bundles on the projective plane and of their \emph{jumping lines} has a long tradition, dating back to the works of Barth, Dolgachev and Kapranov, and Hulek (see \cite{Barth1977,DK93Arr,DK93JLSC,Hulek1979}). The behaviour of such bundles on lines encodes subtle properties of the underlying configurations of points and curves, and appears naturally in the modern theory of unexpected curves (see e.g. \cite{CHMN2018,FV2019,MTG23}). 

In the present paper we develop a new geometric approach, based on \emph{interpolation matrices}, to describe loci of jumping lines for vector bundles arising from configurations of points in $\PP^2$. Our motivation comes from the interpretation of unexpected curves and hypersurfaces proposed in \cite{DMTG25}, where the determinant of a suitable interpolation matrix captures the locus of points at which the corresponding linear system acquires higher than expected dimension. Here we reinterpret this construction in the dual projective plane, showing that, under suitable numerical assumptions, the same determinant encodes the jumping locus of the bundle associated with the configuration.  

More precisely, given a finite set $Z\subset \PP^2$ of $|Z|=n$ points in general position, we associate to it the arrangement of $n$ lines $A_Z$ in the dual plane and the logarithmic (or syzygy) bundle $\cE$ attached to $A_Z$. for a line $L\subset \PP^2$ we denote by $P_L \in (\PP^2)^\vee$ its dual. We are interested in the locus of lines, where the splitting of $\cE|_L$
is non-balanced. We show that this locus can be recovered as the vanishing of the determinant of an explicit interpolation matrix $M(d,m;Z;B)$ depending on the parameters $(d,m)$ determined by $|Z|$, i.e.
$$|Z|=\begin{cases}
2k+1, &(d,m)=(k,k-1),\\
2k, &(d,m)=(2k-1,2k-1).
\end{cases}$$
In both cases the determinant $\det M(d,m;Z;B)$ defines a curve in the dual plane whose degree agrees with the expected degree of the locus of jumping lines of the first kind, i.e. the lines where the splitting type of $\cE$ deviates minimally from the balanced one (see Section \ref{sec:preliminaries}).  For points in general position this curve is irreducible; for special configurations it becomes reducible, revealing geometric relations among the points.

The paper is organised as follows.  
In Section \ref{sec:preliminaries} we recall the basic notation and the definition of the interpolation 
matrix. Section \ref{sec:jumpingLines} contains the main algebraic construction and the proof that $\det M(d,m;Z;B)=0$ indeed describes the locus of jumping lines of the first kind. In Section \ref{sec:pointsGenereal} we discuss the geometric meaning of the parameters $(d,m)$ and the necessary parity assumptions on $|Z|$, clarifying how the cases $c_1(\cE)=0$ and $c_1(\cE)=-1$ lead respectively to the pairs $(k,k-1)$ and $(2k-1,2k-1)$. Finally, in Section \ref{sec:special} we illustrate the theory by a series of explicit examples, showing how components of various multiplicities appear in the determinant, and how the residual part corresponds to genuinely new curves in the dual plane.

This approach bridges classical bundle theory with the modern language of unexpected hypersurfaces: the interpolation matrix plays the same role as the Jacobian syzygies in the construction of \cite{MTG23}, while its determinant geometrically detects the variation of the splitting type. This framework provides a new tool for understanding degeneracy loci of logarithmic bundles and for detecting new geometric features encoded in the associated interpolation matrices.

\section{Preliminaries}\label{sec:preliminaries}

Let $S=\C[x,y,z]$ be the standard graded polynomial ring and coordinate ring of the projective plane $\PP^2=\PP^2_\C$.
For a graded $S$--module $N$ we denote by $\widetilde N$ the associated coherent sheaf on $\PP^2$. Let $L\simeq \PP^1 \subseteq \PP^2$ be a line. Whenever the restriction $\widetilde N|_L$ is locally free, we may apply
the Grothendieck splitting theorem and write
$$\widetilde N|_L \cong \cO_L(a_1)\oplus\dots\oplus \cO_L(a_r),$$
and call the non-decreasing tuple $(a_1,\dots,a_r)$ the \emph{splitting type} of $\widetilde N$ along $L$.

Let $f\in S$ be a reduced homogeneous polynomial of degree $d\ge2$ and $C=\{f=0\}\subset\PP^2$ the corresponding plane curve. We denote by $f_x,f_y,f_z$ the partial derivatives and by
$$ J_f=(f_x,f_y,f_z)\subset S$$
the Jacobian ideal. The module of Jacobian syzygies,
$$AR(f)=\{(a,b,c)\in S^{\oplus 3}: a f_x+b f_y+c f_z=0\},$$
plays a central role in the study of logarithmic vector fields. For a line arrangement $\mathcal A$ defined by $f=\prod_i \ell_i$, we define
$$D(\mathcal A)=\{\theta\in \Der_\C(S)\mid \theta(f)\in (f)\},$$
the $\C$-module of derivations of $S$ tangent to $\mathcal A$. Here $\Der_\C(S)$ denotes the set of $\C$-linear derivations $\theta\colon S\to S$ satisfying the Leibniz rule $\theta(fg)=\theta(f)g+f\theta(g)$, for $ f,g\in S.$ The sheafification of $D(\cA)$ on $\PP^2$ is denoted by $T\langle \cA\rangle$ and called the \emph{logarithmic bundle} (or sheaf of logarithmic vector fields) along $\cA$. For a general reduced curve $C=\{f=0\}\subset\PP^2$, we define the sheaf $T\langle C\rangle$ as the kernel of the map
$$0\longrightarrow T\langle C\rangle\longrightarrow T_{\PP^2}\xrightarrow{\nabla f}\cO_{\PP^2}(d),$$
that is, the subsheaf of tangent vector fields whose contraction with the gradient of $f$ vanishes on $C$.

Sheafifying the standard presentation of the Jacobian ideal
$$0\to AR(f)\to S(-d+1)^{\oplus 3} \xrightarrow{(f_x,f_y,f_z)} J_f \to 0,$$
and twisting by $\cO_{\PP^2}(d-1)$, we obtain an exact sequence
$$0 \longrightarrow \widetilde{AR(f)}(1-d) \longrightarrow\cO_{\PP^2}^{\oplus 3} \xrightarrow{(f_x,f_y,f_z)}\widetilde{J_f}(d-1) \longrightarrow 0.$$
Comparing this with the defining sequence of $T\langle C\rangle$ and using the Euler sequence on $\PP^2$, we conclude that (up to twist) the middle cohomology of the above complex is precisely $T\langle C\rangle$. For line arrangements this coincides with the usual logarithmic derivation bundle.

In many situations $T\langle C\rangle$ (or $T\langle\cA\rangle$) admits a \emph{Steiner presentation}
$$0 \longrightarrow \cO_{\PP^2}(-1)^{\oplus \alpha}\xrightarrow{M}\cO_{\PP^2}^{\oplus(\alpha+2)}\longrightarrow T\langle\cA\rangle(1)\longrightarrow 0,$$
with $M$ a matrix of linear forms in $x,y,z$. This description simplifies the computation of splitting types by restricting $M$ to lines (see \cite[Section 3]{DK93Arr}.

Let $\cE$ be a rank--$2$ vector bundle on $\PP^2$. For a line $L\subset\PP^2$ we have 
$$\cE|_L\simeq\cO_L(a)\oplus\cO_L(b)$$
with $a\ge b$ and $a+b=c_1(\cE)$. The ordered pair $(a,b)$ is the splitting type of $\cE$ on $L$. In the sequel we shall always consider the normalized form of the bundle associated with a configuration of points $Z$. Although the first Chern class $c_1(\cE_Z)$ equals the sum of the exponents of the corresponding syzygy bundle, it is convenient to replace $\cE_Z$ by a suitable twist $\cE_Z(N)$ so that $c_1(\cE_Z(N))\in\{0,-1\}$. For a general line $L\subset\PP^2$, the splitting type is constant and is called  the \emph{generic splitting type} $(a_0,b_0)$.  In the normalized (semi)stable case, one has $0\leq a_0-b_0\le 1$.

A line $L$ is a \emph{jumping line} of $\cE$ if its splitting type differs from the generic splitting type $(a_0,b_0)$ occurring for a general $L$. The set of points in the dual plane $(\PP^2)^\vee$ corresponding to such lines is called the \emph{jumping locus} of $\cE$ and is denoted by $\Jump(\cE)\subset(\PP^2)^\vee$. For logarithmic bundles this jumping behaviour encodes subtle geometry of the underlying curve or arrangement, and can be read off from the graded pieces of $AR(f)$. For a rank--$2$ bundle $\cE$, the \emph{jumping order} of $\cE$ along $L$ is defined as 
$$\jmath_{\cE}(L)=a(L)-a_0\ge0,$$
where $(a(L),b(L))$ is the splitting on $L$ and $(a_0,b_0)$ is the generic splitting type. The strata
\begin{equation}
\label{eq:k-jumpingLine}
   V_k(\cE)=\{L : \; \jmath_{\cE}(L)\ge k\} 
\end{equation}
form a natural filtration of $\Jump(\cE)$ (compare \cite{DS20}). When $c_1(\cE)$ is odd, one can further consider \emph{jumping lines of the second kind}, detected by the behaviour of $\cE$ on the first infinitesimal neighbourhood $L^{(1)}$.

\begin{definition}
Let $L\subset\PP^2$ be a line defined by $\ell=0$ with ideal sheaf $\cI_L$.  
The \emph{first infinitesimal neighbourhood} (or the double line) of $L$ is the closed subscheme
$$L^{(1)}\subset\PP^2 \quad\text{defined by }\cI_{L^{(1)}}=\cI_L^2.$$
\end{definition}

Equivalently, $L^{(1)}$ is the first thickening order of $L$ in the normal direction,
and a homogeneous form $F$ vanishes on $L^{(1)}$ iff $F\in(\ell)^2$.

Locus of jumping lines of the second kind coincides with the degeneracy locus of a canonical quadratic form and in classical settings it is related to the Schur quadric (see \cite{DK93JLSC} for more details).

Fix a line $L\subset\PP^2$ given by $\ell=ax+by+cz=0$.
Restricting the Jacobian map gives
$$\cO_L^{\oplus3}\xrightarrow{(f_x,f_y,f_z)|_L}\cO_L(d-1),$$
whose kernel is a rank--$2$ vector bundle on $L$, isomorphic (up to twist) to $T\langle C\rangle|_L$. Hence the splitting type on $L$ is determined by the degrees of a minimal syzygy basis for $(f_x,f_y,f_z)$ after restriction to $L$. More precisely, if these degrees are $d_1\le d_2$, then $\cE|_L\simeq\cO_L(-d_1)\oplus\cO_L(-d_2)$. A line is jumping precisely when $d_1$ is smaller than its value for a generic line.

Let $Z=\{P_1,\dots,P_s\}\subset\PP^2$ be a finite set of distinct points
and $B\in\PP^2$ a general point. For integers $d\ge m\ge1$ consider the linear system
$$\cL(d;mB+Z)=\{F\in S_d  : \; F|_Z=0, \text{ $F$ vanishes to order $\ge m$ at $B$}\}.$$
Choose homogeneous coordinates $B=(a_0:a_1:a_2)$ and form the \emph{interpolation matrix} $M(d,m;Z;B)$ whose first $s$ rows are the evaluations $w(P_i)$ of a fixed basis $w$ of $S_d$, and whose remaining ${m+1\choose2}$ rows contain all partial derivatives of order $m-1$ of $w$ evaluated at $(a_0,a_1,a_2)$. When ${d+2\choose2}=s+{m+1\choose2}$ the matrix is square, and its determinant
$$F_{d,m;Z}(B)=\det M(d,m;Z;B)$$
is a bihomogeneous polynomial in $(x:y:z)$ and in the coordinates $(a_0:a_1:a_2)$ of $B$ whose zero locus with respect to $(a_0:a_1:a_2)$ parametrizes points $B$ where $\cL(d;\,mB+Z)$ is nonempty.
\begin{example}
We consider the square interpolation matrix in the case $d=2$ and $m=2$.
Let
$$Z=\{P_1,P_2,P_3\}\subset\PP^2,\qquad P_i=(p_{i0}:p_{i1}:p_{i2}),$$
and let $B=(a_0:a_1:a_2)\in\PP^2$ be a general point. Fix a basis of $S_2$, for example
$w=\{x^2,y^2,z^2,xy,xz,yz\}.$
Then the interpolation matrix $M(2,2;Z;B)$ is the $6\times6$ matrix
$$
M(2,2;Z;B)=
\begin{pmatrix}
p_{10}^2 & p_{11}^2 & p_{12}^2 & p_{10}p_{11} & p_{10}p_{12} & p_{11}p_{12} \\
p_{20}^2 & p_{21}^2 & p_{22}^2 & p_{20}p_{21} & p_{20}p_{22} & p_{21}p_{22} \\
p_{30}^2 & p_{31}^2 & p_{32}^2 & p_{30}p_{31} & p_{30}p_{32} & p_{31}p_{32} \\
2a_0 & 0     & 0     & a_1 & a_2 & 0 \\
0     & 2a_1 & 0     & a_0 & 0   & a_2 \\
0     & 0     & 2a_2 & 0   & a_0 & a_1
\end{pmatrix}.
$$
The first three rows correspond to the vanishing conditions at the points of $Z$, while the last three rows encode the vanishing of all first-order partial derivatives at $B$, that is, the condition that a conic vanishes to order at least $2$ at $B$. If the points of $Z$ are not collinear, then $\det M=0$ defines three lines, each passing through a pair of points $P_i$. If the points of $Z$ are collinear, then $\det M=0$ defines that line with multiplicity $3$.
\end{example}

The \emph{fat--point locus}
$$V_{d,m}(Z)=\{B\in\PP^2\; : \; \dim\cL(d;mB+Z)\neq 0\}$$ satisfies $V_{d,m}(Z)=\{F_{d,m;Z}=0\}$.
If
$$h_{j,B}=\dim\cL(d;jB+Z)-\max\left\{0,\dim S_d-s-\binom{j+1}{2}\right\},$$
then the multiplicity of $B$ on $\{F_{d,m;Z}=0\}$ is at least $\sum_{j=1}^d h_{j,B}$ (see \cite[Theorem 2.3]{DMTG25} for details). When $B$ lies on a line $L$ containing $|L\cap Z|$ points, Bezout’s theorem implies 
$$h_{j,B}\ge\max\{0,\,j+|L\cap Z|-d-1\}$$
for suitable $j$. If the sum of all $h_{j,B}$ is bigger than the degree of $F$ then we conclude that $F_{d,m;Z}=0$. And it gives the existence of a curve with prescribed multiplicity at a general point (see \cite[Remark 3.3]{DMTG25}).

These interpolation constructions can be viewed as a combinatorial counterpart of the syzygy description of $T\langle C\rangle$. Fixing a general line $L$ with dual point $P_L$, a syzygy among $(f_x,f_y,f_z)$ yields a rational map $L\dashrightarrow\PP^2$ whose image is a curve passing through prescribed points and having prescribed multiplicity at $P_L$ (\cite[Theorem 3.1]{MTG23}). The degree of the syzygy coincides with one of the summands in the splitting of $T\langle C\rangle|_L$.
Conversely, for configurations $Z$ arising as geometric or combinatorial traces of $C$, the vanishing of $F_{d,m;Z}(B)$ identifies the locus of points where a degree $d$ system with a fat base point exists, and the corresponding degrees and multiplicities match the splitting data on the bundle side (see \cite[Theorem 1.5]{CHMN2018}). In this way the jumping locus $\Jump(T\langle C\rangle)$ often appears as an elimination or projection of the incidence variety $\mathbb{F}_2$ (see \cite{OSS1980}) in $\PP^2\times(\PP^2)^\vee$ cut out simultaneously by the line equation, the interpolation determinant, and the syzygy relations.

For later reference we recall the local invariants of isolated singular points. Let $C\subset\PP^2$ be a reduced curve and let $P\in C$ be a singular point. The \emph{Milnor number} and \emph{Tjurina number} of $C$ at $P$ are defined by
$$\mu_P=\dim_\C \frac{\C\{u,v\}}{\left(\frac{\partial f}{\partial u},\frac{\partial f}{\partial v} \right)}\qquad
\tau_P=\dim_\C \frac{\C\{u,v\}}{\left(f,\frac{\partial f}{\partial u},\frac{\partial f}{\partial v} \right)},
$$
where $f$ is a local analytic equation of $C$ at $P$. For \emph{quasi-homogeneous} singularities (in particular for the ordinary multiple points that occur in our setting) one has $\mu_P=\tau_P$. If $P$ is an \emph{ordinary singularity}, and $C$ has $k$ distinct branches at $P$, then
$$\mu_P=(k-1)^2.$$
The \emph{total Tjurina number} of $C$ is the sum
$$\tau(C)=\sum_{P\in\Sing(C)} \tau_P,$$
and it measures the global defect of freeness of the associated logarithmic bundle;
in particular, for curves with only nodes, $\tau(C)=\binom{|Z|}{2}$.


\section{Jumping lines and the interpolation approach}\label{sec:jumpingLines}

Let $\cE$ be the logarithmic or syzygy bundle introduced earlier, and denote by $V_1(\cE)\subset(\PP^2)^\vee$ the locus of its jumping lines of the first kind, as defined in \eqref{eq:k-jumpingLine}. In this section we show that the classical locus $V_1(\cE)$ can be equivalently described in terms of the interpolation matrix introduced in the previous section.

  We briefly recall the relevant constructions from \cite{OSS1980}, which we shall use throughout the proof.
    Let $\PP^{2\vee}$ be the Grassmannian of lines in $\PP^2$ and
    $$\mathbb{F}_2=\{(x,\ell)\in \PP^2\times \PP^{2\vee} \mid x\in \ell\}$$
    the incidence (flag) manifold with projections
    $$p:\mathbb{F}_2\to \PP^2,\qquad q:\mathbb{F}_2\to \PP^{2\vee}.$$
    Fix a finite set $Z\subset \PP^2$ and integers $d\ge1$, $m\ge1$. Write $\cI_Z\subset \cO_{\PP^2}$ for the ideal sheaf of $Z$ and set
    $$\cE_d:=\cO_{\PP^2}(d)\otimes \cI_Z,\qquad\text{ so that }\qquad p^*\cE_d=\cO_{\mathbb{F}_2}(d)\otimes p^*\cI_Z.$$ For $x\in \PP^2$ consider the submanifold
    $$\mathbb{B}(x)=\{(y,\ell)\in \mathbb{F}_2\mid x,y\in \ell\}=q^{-1}\left(G(x)\right),$$
    where 
    $$G(x):=\{\ell\in \PP^{2\vee}\mid x\in\ell\}\simeq \PP^1$$ 
    denotes the submanifold of all lines through $x$.

\begin{proposition}
\label{prop:iMatAndJL}
    With the above notation, assume that $\binom{d+2}{2}=|Z|+\binom{m+1}{2}$.
    Then
    $$L\in V_1(\cE)\quad\Longleftrightarrow\quad F_{d,m;Z}(P_L)=0,$$
    that is, a line $L$ is a jumping line of the first kind if and only if the corresponding dual point 
    $B=P_L$ satisfies
    $$\dim \cL(d;mB+Z)\;>\;\edim \cL(d;mB+Z)=0.$$
\end{proposition}
\begin{proof}
    Fix a monomial basis of $S_d$ and set $R:=\C[a_0,a_1,a_2]$. Consider the \emph{universal interpolation map}
    $$\Phi_Z\colon S_d\otimes_\C R \longrightarrow R^{|Z|+\binom{m+1}{2}},$$
    whose first $|Z|$ components encode evaluations at the fixed points $Z$ (hence have constant entries in $R$), and whose remaining $\binom{m+1}{2}$ components encode all partial derivatives of order $m-1$ evaluated symbolically at $(a_0,a_1,a_2)$ (hence have polynomial entries in $R$). The matrix of $\Phi_Z$ in the chosen bases is the interpolation matrix $M(d,m;Z)(a_0,a_1,a_2)$ with entries in $R$. For each $B=(a_0:a_1:a_2)\in\PP^2$, the specialization
    \begin{equation}
    \label{eq:kerPhiZB}
        \Phi_{Z,B}:=\Phi_Z\otimes_{R}\C_B\; : \; S_d\longrightarrow \C^{|Z|+\binom{m+1}{2}}
    \end{equation}
    recovers the map at $B$. In particular, in the \emph{square} case
    $$N:=\binom{d+2}{2}=|Z|+\binom{m+1}{2}$$
    we have
    $$\ker\Phi_{Z,B} \cong \cL(d;mB+Z),\qquad F_{d,m;Z}(a_0,a_1,a_2):=\det M(d,m;Z)(a_0,a_1,a_2)\in R,$$
    and
    $$F_{d,m;Z}(B)=0 \Longleftrightarrow \dim_\C\ker\Phi_{Z,B}\ge 1\Longleftrightarrow \cL(d;mB+Z)\neq\{0\}.$$
    Thus $F_{d,m;Z}$ is the bihomogeneous equation on $(\PP^2)^\vee$ cutting out the locus of specialization points $B$ for which the interpolation map drops rank.

   Consider the induced maps
    $$f:=q|_{\mathbb{B}(x)}:\mathbb{B}(x)\to G(x),\qquad\sigma:=p|_{\mathbb{B}(x)}:\mathbb{B}(x)\to \PP^2.$$
    Let $s:G(x)\to \mathbb{B}(x)$ be the section $s(\ell)=(x,\ell)$ (the exceptional section over $x$). Denote by $\cI_s\subset \cO_{\mathbb{B}(x)}$ the ideal sheaf of the image $s(G(x))\subset \mathbb{B}(x)$, and by $\cO_{\mathbb{B}(x)}/\cI_s^m$ the structure sheaf of its $(m-1)$-st infinitesimal neighbourhood.  
    Consider the canonical homomorphism of sheaves
    \begin{equation}
    \label{eq:evalMap}
        f^*f_*\left(\sigma^*\cE_d\right) \longrightarrow\sigma^*\cE_d, 
    \end{equation}
    whose restriction to a fibre $\widetilde L:=f^{-1}(\ell)\simeq L$ is
    $$H^0\left(L,\cO_L(d)\otimes \cI_{Z\cap L}\right)\otimes \cO_L \longrightarrow\cO_L(d)\otimes\cI_{Z \cap L}.$$
    Composing \eqref{eq:evalMap} with the natural quotient gives the jet-evaluation morphism
    $$\Phi_{m,x}: f^*f_*(\sigma^*\cE_d)\longrightarrow \big(\cO_{\mathbb{B}(x)}(d)\otimes\sigma^*\cI_Z\big)/\cI_s^m.$$
    Together with the evaluation map
    $\Phi_Z: f^*f_*(\sigma^*\cE_d)\to\cO_{\mathbb{B}(x)}(d)\otimes\sigma^*\cI_Z,$ we obtain the combined interpolation map
    \begin{equation}
    \label{eq:tensors}
    	\mathbf{\Phi}=(\Phi_Z,\Phi_{m,x}):f^*f_*\left(\sigma^*\cE_d\right)\longrightarrow\left(\cO_{\mathbb{B}(x)}(d)\otimes \sigma^*\cI_Z\right)\oplus\left(\left(\cO_{\mathbb{B}(x)}(d)\otimes \sigma^*\cI_Z\right)\big/\cI_s^{m}\right).
    \end{equation}
    
    Define its kernel sheaf
    $$\widetilde{\cK}:=\ker(\mathbf\Phi)\subset f^*f_*(\sigma^*\cE_d).$$
    Let $\ell\in G(x)$ and $\widetilde L=f^{-1}(\ell)\simeq L$ be the corresponding $f$–fibre. Restricting \eqref{eq:tensors} to $\widetilde L$ and using that $s(G(x))$ meets $\widetilde L$ at the point $x\in L$, we get
    $$H^0\left(L,\cO_L(d) \otimes \cI_{Z\cap L}\right)\otimes \cO_L \longrightarrow \left(\cO_L(d)\otimes \cI_{Z\cap L}\right)\oplus \left(\cO_L(d)\otimes \cO_L/\fm_x^{m}\right),$$
    so that the fibre of $\widetilde{\cK}$ over $\widetilde L$ identifies with
    $$H^0\left(L,\cO_L(d) \otimes \cI_{Z\cap L}\otimes \fm_x^{\,m}\right),$$
    i.e., the space of degree $d$ forms on $L$ vanishing on $Z\cap L$ with multiplicity $\ge m$ at $x$. Equivalently, after choosing coordinates, this is the kernel of the numerical interpolation matrix obtained by evaluating at $Z\cap L$ and taking all partial derivatives of order $m-1$ at $x$.
    Restricting the defining sequence of $\cE$ to a line $L\subset\PP^2$ and recalling that the restriction of the canonical evaluation map \eqref{eq:evalMap} to each fibre $\widetilde L=f^{-1}(\ell)$ is the ordinary interpolation map on $L$, we obtain
    $$H^0(L,\widetilde{\cK}|_{\widetilde L})\simeq\ker\Phi_{Z,B},\qquad B=P_L\in \PP^{2\vee}.$$
    Hence the sheaf $\widetilde{\cK}$ on $\mathbb{B}(x)$ collects fibrewise all kernels of the interpolation maps for lines through $x$, in perfect analogy with the sheaf $q^*q_*p^*\cE\to p^*\cE$ from the standard construction. By cohomology and base change, the formation of $f_*(\sigma^*\cE_d)$ and of its kernel commutes with restriction to fibres. Consequently, the rank of $\mathbf{\Phi}$ drops at a point $\ell\in G(x)$ if and only if
    $$H^0(L,\cO_L(d)\otimes \cI_{Z\cap L}\otimes\fm_x^{m})\neq 0,$$ that is, if and only if the linear system $\cL(d;mB+Z)$ is non-empty for $B=P_L$.
    In the square case $\binom{d+2}{2}=|Z|+\binom{m+1}{2}$ this means
    $$F_{d,m;Z}(B)=0 \Longleftrightarrow \dim\cL(d;mB+Z)>0.$$
    Finally, let $\cE$ be the logarithmic (or syzygy) bundle naturally associated with the data $(Z,d,m)$. By the discussion in \cite[Sec. 2.2]{OSS1980}, a line $L\subset\PP^2$ is a \emph{jumping line of the first kind} for $\cE$ precisely when the restriction $\cE|_L$ acquires non-generic splitting type, i.e. when $H^0(L,\cE|_L(t_0))\ne0$ for the minimal $t_0$ such that the generic restriction has no sections. Since $H^0(L,\widetilde{\cK}|_{\widetilde L})\simeq\ker\Phi_{Z,B}$, this occurs exactly when $\Phi_{Z,B}$ drops rank. Thus the degeneracy locus of $\mathbf{\Phi}$ on $(\PP^2)^\vee$ coincides with the classical jumping locus $V_1(\cE)$
    $$L\in V_1(\cE) \quad \Longleftrightarrow \quad B=P_L\text{ satisfies }F_{d,m;Z}(B)=0,$$
    and the curve $\{F_{d,m;Z}=0\}\subset(\PP^2)^\vee$ describes all
    positions of $B$ for which the corresponding line $L$ is jumping.
\end{proof}

\begin{remark}[Parity condition and choice of parameters]
\label{rem:parity}
    In the proof we implicitly use that the interpolation matrix $M(d,m;Z)(B)$ linearises the evaluation map controlling the restriction $\cE|_L$. This determines the admissible pair $(d,m)$ depending only on the parity of $|Z|$.
    
    \begin{itemize}
    \item If $|Z|=2d'+1$ (so $c_1(\cE)=0$ and $\cE|_L\simeq\cO_L\oplus\cO_L$), the correct interpolation matrix is of type $(d',d'-1)$. Its determinant vanishes exactly on the jumping lines of the first kind.
    \item If $|Z|=2d'$ (so $c_1(\cE)=-1$ and $\cE|_L\simeq\cO_L\oplus\cO_L(-1)$), the appropriate choice is the type $(2d'-1,2d'-1)$. In this almost balanced case, the determinant of the corresponding square matrix again describes the jumping locus of the first kind.
    \end{itemize}
    
\end{remark}

\begin{remark}
    In Proposition \ref{prop:iMatAndJL} we restrict to the locus $V_1(\cE)$ of first-kind jumping lines. For Hulek’s second kind condition the corresponding interpolation matrix becomes rectangular, and its ``rank drop locus'' describes a degeneracy scheme rather than a single curve, so the description via a unique determinant equation does not apply directly.
\end{remark}


\section{Points in general position}\label{sec:pointsGenereal}

As explained in Remark \ref{rem:parity}, the relevant interpolation matrices correspond to specific pairs $(d,m)$ determined by the parity of $|Z|$. In order to avoid unnecessary notation, we shall henceforth use the symbols $(d,m)$ exclusively in this restricted sense, with the precise values specified in each context.

\begin{proposition}\label{prop:s=2d+1 and m=d-1}
    Let $d$ be an integer such that $d \geq 1$ and let $Z\subset\PP^2$ be a set of $s=2d+1$ points in general position. Fix $m=d-1$. Consider the interpolation determinant
    $$F_{d,d-1;Z}(B)=\det M(d,d-1;Z;B),\qquad B\in(\PP^2)^\vee.$$
    Then:
    \begin{itemize}
    \item[(i)] $F_{d,d-1;Z}$ is \emph{not} identically zero; equivalently, for a general line $L$ (i.e. a general $B=P_L$) the map
    $$\Phi_{Z,B}:\; S_d\longrightarrow \C^{|Z|+\binom{d}{2}}$$
    has maximal rank and $\cL(d;(d-1)B+Z)=\{0\}$.
    \item[(ii)] The zero-locus $\{F_{d,d-1;Z}=0\}\subset (\PP^2)^\vee$ is a (reduced) curve of degree
    $$\deg\{F_{d,d-1;Z}=0\}=d(d-1).$$
    Equivalently, for the logarithmic/syzygy rank-two bundle $\cE$ naturally attached to $(d,m;Z)$, the classical jumping locus $V_1(\cE)$ coincides with $\{F_{d,d-1;Z}=0\}$ and has degree equal to its second Chern class,  $c_2(\cE)=d(d-1)$.
    \end{itemize}
\end{proposition}
\begin{proof}
    Set $s=2d+1$ and $m=d-1$. Fix a monomial basis of $S_d$ and let
    $$\Phi_Z: S_d\otimes_\C R\longrightarrow R^{s+\binom{d}{2}},\qquad R=\C[a_0,a_1,a_2],$$
    be the universal interpolation map whose specialization at $B=(a_0:a_1:a_2)$ is
    $$\Phi_{Z,B}: S_d \longrightarrow\C^{s+ \binom{d}{2}},$$
    with the first $s$ components given by evaluation at $Z$ (constants in $B$) and the remaining $\binom{d}{2}$ components given by all partial derivatives of order $d-2$ at $B$ (linear in $B$). In our case
    $$\dim S_d=\binom{d+2}{2}=s+\binom{d}{2},$$
    so $\Phi_{Z,B}$ is square for every $B$. For $Z$ in general position, the $s$ evaluation conditions at $Z$ and the $\binom{d}{2}$ partial derivatives conditions at a general point $B$ are independent. Equivalently, the matrix $M(d,d-1;Z;B)$ has full rank on a nonempty Zariski open subset $U\subset(\PP^2)^\vee$ (upper semicontinuity of the rank). Hence the determinant $F_{d,d-1;Z}$ is not identically zero, and for $B\in U$ one has $\ker\Phi_{Z,B}=0$, i.e.
    $$\cL\left(d;(d-1)B+Z\right)=\{0\}.$$
    That is, for a general position of $B$ the system is empty (its only element is the zero polynomial), while the special points $B$ for which $\cL(d;(d-1)B+Z)\neq\{0\}$ form exactly the degeneracy locus $\{F_{d,d-1;Z}=0\}$.
    This proves (i).

    Since $F_{d,d-1;Z}\not\equiv 0$, its zero set $\{F_{d,d-1;Z}=0\}\subset(\PP^2)^\vee$ is a divisor. We compute its degree. Let $\cE$ be the logarithmic/syzygy rank-two bundle associated with the data $(d,m;Z)$. By Theorem \ref{prop:iMatAndJL} the degeneracy curve of the square interpolation map coincides with the jumping locus $V_1(\cE)$, and its degree equals $c_2(\cE)$. For $Z$ consisting of $2d+1$ general (node-type) points one has
    $$c_2(\cE)=3d^2-\tau(Z)=3d^2-\binom{2d+1}{2}=d(d-1).$$
    On the other hand, as shown for the interpolation determinant in
    \cite[Theorem 2.3]{DMTG25} (or equivalently by degree counting in the matrix
    $M(d,d-1;Z;B)$), one has $\deg F_{d,d-1;Z}=d(d-1)$. Therefore the two invariants coincide, $\deg\{F_{d,d-1;Z}=0\}=c_2(\cE)=d(d-1),$ which confirms the statement of the proposition.
\end{proof}
\begin{remark}
    While the existence of a jumping curve of degree $d(d-1)$ for $c_1=0$ 
    fits into the general framework of curves of jumping lines (see \cite[Section 2.2.4]{OSS1980}), we are not aware of any explicit interpolation matrix realization of this curve in the literature. Our determinantal description therefore appears to be new. See also \cite[Remark 6.4]{DS20} for related predictions in the stable case.
\end{remark}

\begin{proposition}\label{prop:s=d+1 and m=d}
    Let $d=2k-1$, for positive integer $k$, and let $Z=\{P_1,\dots,P_{d+1}\}\subset\PP^2$ be the set of points in general position. Fix $m=d$ and consider the square interpolation determinant
    $$F_{d,d;Z}(B)=\det M(d,d;Z;B),\qquad B\in(\PP^2)^\vee.$$
    Then the degeneracy locus is the union of all secant lines to $Z$:
    $$ \{F_{d,d;Z}=0\}=\bigcup_{1\le i<j\le d+1}\langle P_i,P_j\rangle \subset (\PP^2)^\vee.$$
    In particular, $F_{d,d;Z}$ splits (over $\C$) as a product of linear forms defining these $\binom{d+1}{2}$ lines.
\end{proposition}
\begin{proof}
    Let $Q\in\PP^2$. If a degree $d$ curve $C$ has multiplicity $\ge d$ at $Q$, then in affine coordinates $(u,v)$ centered at $Q$ its local equation has no terms of degree $<d$, hence equals its tangent cone (the homogeneous degree $d$ part). In two variables every homogeneous degree $d$ polynomial splits over $\C$ into linear forms, so $C$ is a cone of degree $d$ with vertex $Q$, i.e.
    $$C=L_1+\ldots+L_d\qquad \text{ with each } L_i \text{ a line through } Q.$$
    Assume $\cL(d;\,dQ+Z)\neq\{0\}$ for some $Q\in\PP^2$. Then every $C\in\cL(d;\,dQ+Z)$ is a sum of $d$ lines through $Q$. Since $|Z|=d+1$, to cover all points $Z=\{P_1,\dots,P_{d+1}\}$ with $d$ lines through a single point $Q$, at least one of these lines must contain two points of $Z$. Thus $Q$ lies on a secant line $\langle P_i,P_j\rangle$ for some $i\ne j$. Equivalently, the corresponding point $B=P_{\langle P_i,P_j\rangle}\in(\PP^2)^\vee$ belongs to the degeneracy locus of $F_{d,d;Z}$.

    Conversely, let $Q\in \langle P_i,P_j\rangle$ be any point on a secant line to $Z$.
    Define a cone
    $$C=\langle P_i,P_j\rangle \; +\sum_{r\in\{1,\ldots,d+1\}\setminus\{i,j\}}\langle Q,P_r\rangle.$$
    Then $C$ is of degree $d$, it has multiplicity $d$ at $Q$ (sum of $d$ lines through $Q$), and it passes through every point of $Z$. Hence $C\in\cL(d;\,dQ+Z)$, so the interpolation matrix $M(d,d;Z;B)$ drops rank at the dual point $B=P_{\langle P_i,P_j\rangle}$.

    Therefore, the degeneracy locus consists precisely of the points of $(\PP^2)^\vee$ representing the secant lines $\langle P_i,P_j\rangle$ to $Z$, i.e.
    $$\{F_{d,d;Z}=0\}=\left\{P_{\langle P_i,P_j\rangle}: 1\le i<j \le d+1 \right\} = \bigcup_{1\le i<j \le d+1}\langle P_i,P_j\rangle \subset (\PP^2)^\vee,$$
    where, by abuse of notation, we identify a line in $\PP^2$ with the corresponding point in the dual plane. Since $Z$ is in general position, no three points are collinear; therefore the secants are pairwise distinct and appear with multiplicity one. It follows that $F_{d,d;Z}$ factors (over $\C$) as the product of the $\binom{d+1}{2}$ linear forms defining these points in $(\PP^2)^\vee$, as claimed.
\end{proof}

\begin{remark}
    The description above matches the classical picture for Hulsbergen bundles with $c_1=-1$, $c_2=|Z|=d+1$, where the jumping lines form the reduced union of all secants through the base points (see \cite[Section 10]{Hulek1979} for details). For $c_1=0$, Barth proved in \cite{Barth1977} that a stable rank-two bundle on $\PP^2$ is determined by its curve of jumping lines together with a $\theta$-characteristic, and analyzed which curves can occur for small values of $c_2$ (see the survey in \cite[Section 2.2.4]{OSS1980}). 
\end{remark}

\begin{remark}
When the points $Z$ are in general position, the determinant $M(d,m;Z)(B)$ defines an irreducible curve in the dual plane $\PP^{2\vee}$. Indeed, in this case there are no special linear or quadratic relations among the evaluation conditions, so no fixed components appear in the determinant. Hence the corresponding jumping locus is an irreducible curve of the expected degree.

For special configurations (e.g. when six points of $Z$ lie on a conic and five on two lines, see Example \ref{ex:11pts+conic+lines}) the determinant may factor, and some fixed components can occur, but this does not affect the general construction or its motivation.
\end{remark}


\section{Special configurations of points and forced components of the degeneracy locus}\label{sec:special}

Throughout this section let $C\subset\PP^2$ be an irreducible curve of degree $t\ge 2$, and set
$$|C|:=|Z\cap C|,$$
for a number of points on $C$ from the fixed set of points $Z$. Following \cite[Theorem 2.3]{DMTG25} for a point $B\in\PP^2$ and an integer $j\ge 1$ we write
$$h_{j,B}:=\dim \cL(d;jB+Z)-\vdim \cL(d;jB+Z),$$
so that $h_{j,B}$ measures the decrease of the interpolation rank in order $j$ at $B$. We will bound $h_{j,B}$ uniformly for all $B\in C$ under natural “fixed component” hypotheses for $C$. The following proposition is a generalization of result and method that is provided in \cite[Remark 3.3]{DMTG25}.

\begin{proposition}\label{prop:degC hjB bound}
    Let $d\ge2$, $d \ge m\ge t$, and let $C\subset\PP^2$ be an irreducible curve of degree $t\ge 1$. Assume that for a given integer $j\ge t$, every curve in $\cL(d; jB+Z)$ contains $C$ as a fixed component for all $B\in C$. 
    Then, for every $B\in C$ and for each integer $j$ satisfying
    $$j_{\min}:=\max\left\{t,\left\lceil td -|C| -\frac{t^2-3t}{2}+1\right\rceil\right\}\le j\le m,$$
    one has the uniform lower bound
    $$h_{j,B}\ge  |C|+j -td +\frac{t^2-3t}{2}.$$
\end{proposition}
\begin{proof}
    Fix $j\ge t$ and $B\in C$. By the hypothesis, every curve in $\cL(d;\,jB+Z)$ has $C$ as a fixed component, hence
    $$\cL(d;jB+Z)\subset (C)+\cL\left(d-t;(j-1)B+Z\setminus C\right).$$
    Consequently
    $$h_{j,B}\ge\vdim \cL\left(d-t;(j-1)B+Z\setminus C\right)-\vdim \cL\left(d;jB+Z\right).$$
    Using $\vdim \cL(d;jB+Z)=\binom{d+2}{2}-|Z|-\binom{j+1}{2}$ and similarly for the $(d-t)$–system, we get
    \begin{align*}
    h_{j,B} &\ge\left(\binom{d-t+2}{2}-|Z\setminus C|-\binom{j}{2}\right) -\left(\binom{d+2}{2}-|Z|-\binom{j+1}{2}\right)\\
    &=\frac{-t(2d+3-t)}{2}+j+|C| =|C|+j -td +\frac{t^2-3t}{2}.
    \end{align*}
    Finally, the inequality $|C|+j -td +\frac{t^2-3t}{2}\ge 1$ is equivalent to 
    $j\ge td -|C| -\frac{t^2-3t}{2}+1, $
    from which, after rounding up, the statement follows.
\end{proof}
\begin{remark}
    Although Proposition \ref{prop:degC hjB bound} can be formulated and proved for any $t\ge 1$, in the geometric situation relevant for the description of the jumping--line locus we always assume $s=|Z|=2d+1$. Under this assumption the range of admissible integers $j$ in the proposition is non--empty only for $t=1,2$.
    
    Indeed, for $t\ge 3$ the inequalities in the statement force
    $$d-1=m \ge j \ge j_{\min}\ge td-|C|-\frac{t^{2}-3t}{2}+1.$$
    Ignoring $j$ and requiring only that $m\ge j_{\min}$ gives
    \begin{equation}
    \label{eq:lowerBound}
        |C|\ge (t-1)d-\frac{t^{2}-3t}{2}+2.
    \end{equation}
    Since $C\subset Z$ we have $|C|\le s=2d+1$, hence
    $$2d+1 \ge(t-1)d-\frac{t^{2}-3t}{2}+2,$$
    or equivalently
    $t^{2}-3t-2 \ge 2(t-3)d.$
    
    On the other hand, if a curve of degree $t$ is to occur as a fixed component of a configuration of $2d+1$ points, Bézout’s theorem imposes the necessary condition
    $$2d+1=s\ge\binom{t+2}{2}, \qquad \text{that is,} \qquad d\ge\frac{t^{2}+3t}{4}.$$
    Substituting this into the previous inequality yields
    $-t^{3}+2t^{2}+3t-4\ge 0,$
    which holds only for $t=1,2$. Thus no curves of degree $t\ge 3$ can appear as fixed components in the situation of Proposition \ref{prop:degC hjB bound} under the assumptions relevant for the jumping--line locus.
\end{remark}
\begin{example}\label{ex:11pts+conic+lines}
    Let $Z\subset\PP^2$ consists of nine points arranged as follows: seven of them lie on a smooth conic $C$, while the remaining two points on line $L$, which also contains two points from conic. The configuration is illustrated in the adjacent figure.

    \vspace{0.2cm}

    \noindent
    \begin{minipage}{0.55\textwidth}
    Since $|Z|=9=2\cdot4+1$, by the parity rule of Remark \ref{rem:parity} we consider the interpolation matrix of type $(4,3)$. According to Proposition \ref{prop:iMatAndJL}, the determinant $F=\det M(4,3;Z;B)$ defines the jumping locus in the dual plane, provided it is not identically zero. The expected degree is
    $$\deg(F)=2\binom{4}{2}=12.$$

    \end{minipage}
    \begin{minipage}{0.05\textwidth}
    \hfill
    \end{minipage}
    \begin{minipage}{0.30\textwidth}
    \begin{tikzpicture}[line cap=round,line join=round,>=triangle 45,x=1.0cm,y=1.0cm]
    \clip(-0.12,-2.58) rectangle (5.72,2.54);
    \draw [rotate around={-175.4503851258514:(2.8296569803254066,-0.17472291176206367)},line width=0.5pt,dash pattern=on 2pt off 2pt] (2.8296569803254066,-0.17472291176206367) ellipse (1.9189627118883537cm and 1.6314344340769003cm);
    \draw [line width=0.5pt,dash pattern=on 2pt off 2pt,domain=-0.12:5.72] plot(\x,{(-3.4319256682486894--0.8246717381473854*\x)/3.7239073090466834});
    \begin{scriptsize}
    \draw [fill=black] (0.94,0.06) circle (1.5pt);
    \draw [fill=black] (2.22,1.36) circle (1.5pt);
    \draw [fill=black] (4.38,0.82) circle (1.5pt);
    \draw [fill=black] (4.56,-0.84) circle (1.5pt);
    \draw [fill=black] (2.04,-1.68) circle (1.5pt);
    \draw [fill=black] (4.722960246344231,0.12432268820985493) circle (1.5pt);
    \draw [fill=black] (0.9990529372975474,-0.7003490499375304) circle (1.5pt);
    \draw [fill=black] (2.3216509972837023,-0.40745528259014646) circle (1.5pt);
    \draw [fill=black] (3.575456524994442,-0.12979585177560415) circle (1.5pt);
    \end{scriptsize}
    \end{tikzpicture}
    \end{minipage}

    \vspace{0.2cm}

    We now analyse the possible fixed components of $F$ using Proposition \ref{prop:degC hjB bound}. For the conic $C$ we have
    $$j_{\min}=\max\left\{2,\left\lceil 8-7+1+1 \right\rceil\right\}=3,$$
    so for each point $B\in C$,
    $$h_{3,B}\ge 7+3-8-1=1.$$
    Hence $C$ appears in the divisor of $F$ with multiplicity at least $1$.
    
    For a point $B\in L$ we have 
    $$j_{\min}=\max\left\{1,\left\lceil 8-7+1+1 \right\rceil\right\}=3,$$
    and consequently
    \begin{equation*}
    h_{2,B}\ge 4+2-4-1=1,\qquad h_{3,B}\ge 4+3-4-1=2,
    \end{equation*}
    so the line $L$ occurs in the locus with multiplicity at least $1+2=3$.
    
    Altogether, if the determinant $F$ is not identically zero, the jumping locus contains
    $$C+3L$$
    as fixed components, contributing a total degree of $5$. Therefore, the residual part of $F$ must have degree $7$, corresponding to an additional, genuinely new component of the jumping locus. A direct computation for a specific configuration of nine points with the above incidence type shows that
    $$F=F_C \cdot F_L^3 \cdot G,$$
    where $F_C$ and $F_L$ define $C$ and $L$, and $G$ is a form of degree $7$ which is not divisible by $F_C$ nor by $F_L$. In particular, for this configuration the jumping locus is the union of $C$, the triple line $3L$, and an additional component of degree $7$. 
    
    We have performed the computation in {\sc Singular} \cite{singular}; the corresponding script 
    is available in \cite{JLgithub}. Moreover, running the same script on many randomly chosen configurations with the same incidence type always produced the same factorization pattern, which provides strong experimental evidence that this behaviour is generic.
\end{example}
\begin{remark}
    In the setting of Proposition \ref{prop:degC hjB bound} we cannot obtain, even in principle, a configuration whose jumping locus contains two distinct conics as fixed components. Indeed, for $t=2$ the lower bound in \eqref{eq:lowerBound} from the previous remark reduces to $|C| \ge d+3.$ Assume now that two irreducible conics $C_{1}$ and $C_{2}$ occur as components of the jumping locus. Since $|Z|=2d+1$ and two conics meet in
    $4$ points (and fewer intersection points only worsen the inequality), we have
    $|C_{1}|+|C_{2}|=2d+1.$ The bound \eqref{eq:lowerBound} forces $|C_{1}|\ge d+3$ and
    $|C_{2}|\ge d+3$, which is impossible since $2d+6>2d+1$. Even in the minimal case $|C_{1}|=d+3$ and $|C_{2}|=d+2$ the latter conic fails to satisfy \eqref{eq:lowerBound}, and therefore cannot be a fixed component of the degeneracy locus.
    
    Thus, under the hypotheses of Proposition \ref{prop:degC hjB bound}, at most one conic may occur as a forced component of the jumping locus.
\end{remark}
\begin{remark}
    Computer experiments strongly suggest that the above restriction is not merely an artefact of Proposition \ref{prop:degC hjB bound}. For configurations in which the points are distributed on two conics with $|C_{1}|=d+3$ and $|C_{2}|=d+2$, explicit computations of the interpolation matrices $M(d,m;Z)$ consistently show that only one of the two conics appears in the factorization of $\det M(d,m;Z)$. This pattern persists for all randomly generated configurations of the same incidence type tested by the authors.
    
    The scripts used for these computations (implemented in {\sc Singular} \cite{singular}) are available for independent verification in \cite{JLgithub}.
\end{remark}


\subsection*{Acknowledgements.}
The authors are grateful to Tomasz Szemberg for insightful remarks and constructive feedback on a preliminary version of this article, and to Piotr Pokora for reading the manuscript.

\noindent
Guardo and Keiper's research is partially supported by Gnsaga of Indam, Progetto Piaceri 2024-26, Università di Catania and by PRIN 2022, “$0$-dimensional schemes, Tensor Theory and applications” – funded by the European Union Next Generation EU, Mission 4, Component 2 -- CUP: E53D23005670006.

\noindent
Grzegorz Malara was supported by the National Science Centre (Poland), Sonata Grant \\ \textbf{2023/51/D/ST1/00118}. A significant part of this work was carried out during the author’s visit to the Università di Catania, which the author gratefully acknowledges for its hospitality.

 
\printbibliography


\footnotesize
\noindent
Elena Guardo, Graham Keiper: \textsc{Dipartimento di Matematica e Informatica, Viale A. Doria 6, 95125, Catania, Italy}\\
\textit{E-mail address:} \texttt{guardo@dmi.unict.it}\\
\textit{E-mail address:} \texttt{gtrkeiper@outlook.com}\\
	
\noindent
Grzegorz Malara: \textsc{Department of Mathematics, University of the National Education Commission, Podchor\c a\.zych 2, 30-084 Krak\'ow, Poland}\\
\textit{E-mail address:}
\texttt{grzegorzmalara@gmail.com}\\

\end{document}